\documentclass[12pt]{article}
\usepackage{amsmath}
\usepackage{amsfonts}
\usepackage{mitpress}
\usepackage{graphicx}

\begin{document}

\title{Galerkin Method for the numerical solution of the RLW equation by
using exponential B-splines}
\author{M. Z. G\"{o}rg\"{u}l\"{u}, \.{I}. Da\u{g} and D. Irk \\
{\small Department of Mathematics-Computer Science, Eskisehir Osmangazi
University, 26480, Eskisehir, Turkey.}}
\maketitle

\begin{abstract}
In this paper, the exponential B-spline functions are used for the numerical
solution of the RLW equation. Three numerical examples\ related to
propagation of single solitary wave, interaction of two solitary waves and
wave generation are employed to illustrate the accuracy and the efficiency
of the method. Obtained results are compared with some early studies.

\end{abstract}

\section{Introduction}

Many studies exist for the numerical solutions of the differential equations
using splines. Splines are piecewise functions which have certain continuity
at the joint points given to set up the splines. The spline related
numerical techniques mainly offer the economical computer code and easy
computational calculations. Thus they are preferable in forming the
numerical methods. Until now, polynomial splines have been extensively
developed and used for approximation of curve and surfaces and finding
solutions of the differential equations. The polynomial spline based
algorithms have been found to be quite advantageous for finding solutions of
the differential equations. Because it has been demonstrated that they yield
the lower cost and simplicity to write the program code. Base of splines
known as the B-splines is also widely used to build up the trial functions
for numerical methods. The exponential spline is proposed to be more general
form of these splines. In the approximation theory, the exponential
B-splines are shown to model the data which have\ sudden growth and decay\
whereas polynomials are not appropriate due to having osculatory behavior.
Since some differential equations have steep solutions, the use of the
exponential B-splines in the numerical methods may exhibit good solutions
for differential equations. McCartin \cite{mccartin} has introduced the
exponential B-spline as a basis for the space of exponential splines. The
exponential B-spline properties accord with those of polynomial B-splines
such as smootness, compact support, positivity, recursion for derivatives.
Thus the exponential B-splines can be used as the trial function for the
variational methods such as Galerkin and collocation methods.

The exponential B-spline based methods have been started to solve some
differential equations: Numerical solution of the singular perturbation
problem is solved with a variant of exponential B-spline collocation method
in the work \cite{ms1}, the cardinal exponential B-splines is used for
solving the singularly perturbed problems \cite{ca}, the exponential
B-spline collocation method is built up for finding the numerical solutions
of the self-adjoint singularly perturbed boundary value problems in the work
\cite{rao}, the numerical solutions of the Convection-Diffusion equation is
obtained by using the exponential B-spline collocation method \cite%
{mohammadi}.

The collocation methods based on the exponential B-spline functions have
been constructed to solve the differential equations. In this study, the
exponential B-spline function are used to set up the trial functions which
are placed in place of the unknown variable of the differential equations \
for the Galerkin finite element method. Thus nonlinear RLW equation will be
solved with the proposed method numerically. The RLW equation describes a
large number of important physical phenomena, such as shallow waters and
plasma waves. Therefore it plays a major role in the study of nonlinear
dispersive waves. Because of having limited analytical solutions, numerical
analysis of the RLW equation has an importance in its study.

Various techniques have been developed to obtain the numerical solution of
this nonlinear partial differential equation, some of which are finite
difference methods \cite{peregrine,iskandar,jain,shankar,esen}, finite
element methods \cite%
{alexander,serna,luo,saka2,saka3,dogan,dag2,dag3,saka,avilez,mei} and
spectral methods \cite{guo,sloan,manoranjon}.

The paper is outlined as follows. In section 2, exponential B-splines and
their some basic relations are introduced. In section 3, the application of
the numerical method is given. The efficiency and the accuracy of the
present method are investigated by using three numerical experiments related
to propagation of single solitary wave, interaction of two solitary waves
and wave generation. Finally some remarks are concluded in the last section.

\section{Exponential B-splines and Finite Element Solution}

In this study, we will consider the regularized long wave (RLW) equation
\begin{equation}
u_{t}+u_{x}+\varepsilon uu_{x}-\mu u_{xxt}=0  \label{1}
\end{equation}%
where $x$ is space coordinate, $t$ is time, $u$ is the wave amplitude and $%
\varepsilon $ and $\mu $ are positive parameters. Boundary and initial
conditions of the Eq.(\ref{1}) are

\begin{equation}
\begin{tabular}{l}
$u\left( a,t\right) =\beta _{1},u\left( b,t\right) =\beta _{2}$ \\
$u_{x}\left( a,t\right) =0,\text{ }u_{x}\left( b,t\right) =0,t\in \left( 0,T%
\right] ,$ \\
$u\left( x,0\right) =f\left( x\right) ,x\in \left[ a,b\right] .$%
\end{tabular}
\label{2}
\end{equation}

This equation (\ref{1}) was first introduced by Peregrine \cite{peregrine}
for modelling the propagation of unidirectional weakly nonlinear and weakly
dispersive water waves. Since the RLW equation obviates the certain
problematical aspects of KdV equation and it generally has more expedient
mathematical properties, Benjamin et al. \cite{benjamin} proposed the use of
RLW equation as a preferable model to the KdV equation.

Let us consider a uniform mesh $\Gamma $ with the knots $x_{i}$ on $\left[
a,b\right] $ such that%
\begin{equation*}
\Gamma :a=x_{0}<x_{1}<x_{2}<\cdots <x_{N-1}<x_{N}=b
\end{equation*}%
\noindent where $h=\dfrac{b-a}{N}$ and $x_{i}=x_{0}+ih.$

Let $B_{i}\left( x\right) $ be the B-splines at the points of $\Gamma $
together with knots $x_{i}$ , $i=-3,-2,-1,N+1,N+2,N+3$ outside the interval $%
[a,b]$ and having a finite support on the four consecutive intervals $\left[
x_{i}+kh,x_{i}+\left( k+1\right) h\right] _{k=-3}^{0},$ $i=0,...,N+2.$ The $%
B_{i}\left( x\right) $ can be defined as

\begin{equation}
B_{i}\left( x\right) =\left\{
\begin{array}{lll}
b_{2}\left[ \left( x_{i-2}-x\right) -\dfrac{1}{p}\left( \sinh \left( p\left(
x_{i-2}-x\right) \right) \right) \right] & \text{ \ } & \text{if }x\in \left[
x_{i-2},x_{i-1}\right] ; \\
a_{1}+b_{1}\left( x_{i}-x\right) +c_{1}e^{p\left( x_{i}-x\right)
}+d_{1}e^{-p\left( x_{i}-x\right) } & \text{ } & \text{if }x\in \left[
x_{i-1},x_{i}\right] ; \\
a_{1}+b_{1}\left( x-x_{i}\right) +c_{1}e^{p\left( x-x_{i}\right)
}+d_{1}e^{-p\left( x-x_{i}\right) } & \text{ } & \text{if }x\in \left[
x_{i},x_{i+1}\right] ; \\
b_{2}\left[ \left( x-x_{i+2}\right) -\dfrac{1}{p}\left( \sinh \left( p\left(
x-x_{i+2}\right) \right) \right) \right] & \text{ } & \text{if }x\in \left[
x_{i+1},x_{i+2}\right] ; \\
0 & \text{ } & \text{otherwise.}%
\end{array}%
\right.  \label{3}
\end{equation}%
where

$%
\begin{array}{l}
p=\underset{0\leq i\leq N}{\max }p_{i},\text{ }s=\sinh \left( ph\right) ,%
\text{ }c=\cosh \left( ph\right) \\
b_{2}=\dfrac{p}{2\left( phc-s\right) },\text{ }a_{1}=\dfrac{phc}{phc-s},%
\text{ }b_{1}=\dfrac{p}{2}\left[ \dfrac{c\left( c-1\right) +s^{2}}{\left(
phc-s\right) \left( 1-c\right) }\right] , \\
c_{1}=\dfrac{1}{4}\left[ \dfrac{e^{-ph}\left( 1-c\right) +s\left(
e^{-ph}-1\right) }{\left( phc-s\right) \left( 1-c\right) }\right] ,\text{ }%
d_{1}=\dfrac{1}{4}\left[ \dfrac{e^{ph}\left( c-1\right) +s\left(
e^{ph}-1\right) }{\left( phc-s\right) \left( 1-c\right) }\right]%
\end{array}%
$

Each basis function $B_{i}\left( x\right) $ is twice continuously
differentiable. The values of $B_{i}\left( x\right) ,$ $B_{i}^{\prime
}\left( x\right) $ and $B_{i}^{\prime \prime }\left( x\right) $ at the knots
$x_{i}$'s are obtained from the Table 1.

\begin{equation*}
\begin{tabular}{c}
{\small Table 1: Exponential B-spline values} \\ \hline
\multicolumn{1}{|c|}{%
\begin{tabular}{lccccc}
& $x_{i-2}$ & $x_{i-1}$ & $x_{i}$ & $x_{i+1}$ & $x_{i+2}$ \\ \hline
$B_{i}\left( x\right) $ & $0$ & $\frac{s-ph}{2\left( phc-s\right) }$ & $1$ &
$\frac{s-ph}{2\left( phc-s\right) }$ & $0$ \\
$B_{i}^{\prime }\left( x\right) $ & $0$ & $\frac{p\left( 1-c\right) }{%
2\left( phc-s\right) }$ & $0$ & $\frac{p\left( c-1\right) }{2\left(
phc-s\right) }$ & $0$ \\
$B_{i}^{\prime \prime }\left( x\right) $ & $0$ & $\frac{p^{2}s}{2\left(
phc-s\right) }$ & $\frac{-p^{2}s}{phc-s}$ & $\frac{p^{2}s}{2\left(
phc-s\right) }$ & $0$%
\end{tabular}%
} \\ \hline
\end{tabular}%
\end{equation*}

The $B_{i}\left( x\right) ,i=-1,\ldots ,N+1$ form a basis for functions
defined on the interval $[a,b]$. We seek an approximation $U_{N}(x,t)$ to
the analytical solution $u(x,t)$ in terms for the exponential B-splines
\begin{equation}
u\left( x,t\right) \approx U_{N}\left( x,t\right) =\overset{N+1}{\underset{%
i=-1}{\sum }}B_{i}\left( x\right) \delta _{i}\left( t\right)  \label{4}
\end{equation}%
where $\delta _{i}\left( t\right) $ are time dependent unknown to be
determined from the boundary conditions and Galerkin approach to the
equation (\ref{1}). The approximate solution and their derivatives at the
knots can be found from the Eq. (\ref{3}-\ref{4}) as
\begin{equation}
\begin{tabular}{l}
$U_{i}=U_{N}^{{}}(x_{i},t)=\alpha _{1}\delta _{i-1}+\delta _{i}+\alpha
_{1}\delta _{i+1},$ \\
$U_{i}^{\prime }=U_{N}^{\prime }(x_{i},t)=\alpha _{2}\delta _{i-1}-\alpha
_{2}\delta _{i+1},$ \\
$U_{i}^{\prime \prime }=U_{N}^{\prime \prime }(x_{i},t)=\alpha _{3}\delta
_{i-1}-2\alpha _{3}\delta _{i}+\alpha _{3}\delta _{i+1}$%
\end{tabular}
\label{5}
\end{equation}%
where $\alpha _{1}=\dfrac{s-ph}{2(phc-s)},\alpha _{2}=\dfrac{p(1-c)}{2(phc-s)%
},\alpha _{3}=\dfrac{p^{2}s}{2(phc-s)}.$

Applying the Galerkin method to the RLW equation with the exponential
B-splines as weight function over the element $\left[ a,b\right] $ gives
\begin{equation}
\underset{a}{\overset{b}{\int }}B_{i}\left( x\right) \left(
u_{t}+u_{x}+\varepsilon uu_{x}-\mu u_{xxt}\right) dx=0.  \label{6}
\end{equation}

The approximate solution $U_{N}$ \ over the element $[x_{m},x_{m+1}]$ can be
written as

\begin{equation}
U_{N}^{e}=B_{m-1}\left( x\right) \delta _{m-1}\left( t\right) +B_{m}\left(
x\right) \delta _{m}\left( t\right) +B_{m+1}\left( x\right) \delta
_{m+1}\left( t\right) +B_{m+2}\left( x\right) \delta _{m+2}\left( t\right)
\label{7}
\end{equation}%
where quantities $\delta _{j}\left( t\right) ,j=m-1,...,m+2$ are element
parameters and $B_{j}\left( x\right) ,j=m-1,...,m+2$ are known as the
element shape functions.

The contribution of the integral equation (\ref{6}) over the sample interval
$[x_{m},x_{m+1}]$ is given by
\begin{equation}
\underset{x_{m}}{\overset{x_{m+1}}{\int }}B_{j}\left( x\right) \left(
u_{t}+u_{x}+\varepsilon uu_{x}-\mu u_{xxt}\right) dx=0.  \label{8}
\end{equation}%
Applying the Galerkin discretization scheme by replacing approximate
solution $U_{N}^{e}$ (\ref{7}) and its derivatives $(U_{N})_{t},$ $%
(U_{N})_{x},$ $(U_{N})_{xxt}$ into the exact solution $u$ and its
derivatives $u_{t},$ $u_{x},$ $u_{xxt}$ respectively, we obtain a system of
equations in the unknown parameters $\delta _{j}$
\begin{equation}
\begin{array}{l}
\overset{m+2}{\underset{i=m-1}{\sum }}\left\{ \left( \underset{x_{m}}{%
\overset{x_{m+1}}{\int }}B_{j}B_{i}dx\right) \overset{\mathbf{\bullet }}{%
\delta }_{i}+\left( \underset{x_{m}}{\overset{x_{m+1}}{\int }}%
B_{j}B_{i}^{\prime }dx\right) \delta _{i}\right. \\
\left. +\varepsilon \left( \underset{x_{m}}{\overset{x_{m+1}}{\int }}%
B_{j}\left( \overset{m+2}{\underset{k=m-1}{\sum }}\delta _{k}B_{k}\right)
B_{i}^{\prime }dx\right) \delta _{i}-\mu \left( \underset{x_{m}}{\overset{%
x_{m+1}}{\int }}B_{j}B_{i}^{\prime \prime }dx\right) \overset{\mathbf{%
\bullet }}{\delta }_{i}\right\} =0%
\end{array}
\label{9}
\end{equation}%
where $i,j$ and $k$ take only the values $m-1,m,m+1,m+2$ for $m=0,1,\ldots
,N-1$ and $\overset{\mathbf{\bullet }}{}$ denotes time derivative.

In the above system of differential equations, when $A^{e},B^{e},C^{e}(%
\delta )$ and $D^{e}$ are denoted by

\begin{equation}
\begin{tabular}{ll}
$A^{e}=\underset{x_{m}}{\overset{x_{m+1}}{\int }}B_{j}B_{i}dx,$ & $B^{e}=%
\underset{x_{m}}{\overset{x_{m+1}}{\int }}B_{j}B_{i}^{\prime }dx,$ \\
$C^{e}\left( \delta \right) =\underset{x_{m}}{\overset{x_{m+1}}{\int }}%
B_{j}\left( \overset{m+2}{\underset{k=m-1}{\sum }}\delta _{k}B_{k}\right)
B_{i}^{\prime }dx,$ & $D^{e}=\underset{x_{m}}{\overset{x_{m+1}}{\int }}%
B_{j}B_{i}^{\prime \prime }dx$%
\end{tabular}
\label{10}
\end{equation}%
where $A^{e},B^{e}$ and $D^{e}$ are the element matrices of which dimensions
are $4\times 4$ and $C^{e}\left( \delta \right) $ is the element matrix with
the dimension $4\times 4\times 4$, the matrix form of the Eq.(\ref{9}) can
be written as%
\begin{equation}
\left( \mathbf{A}^{e}-\mu \mathbf{D}^{e}\right) \overset{\mathbf{\bullet }}{%
\mathbf{\delta }^{e}}+\left( \mathbf{B}^{e}+\varepsilon \mathbf{C}^{e}\left(
\mathbf{\delta }^{e}\right) \right) \mathbf{\delta }^{e}=0  \label{11}
\end{equation}%
where $\mathbf{\delta }^{e}\mathbf{=}\left( \delta _{m-1},...,\delta
_{m+2}\right) ^{T}$

Gathering the systems (\ref{11}) over all elements, we obtain global system

\begin{equation}
\left( \mathbf{A}-\mu \mathbf{D}\right) \overset{\mathbf{\bullet }}{\mathbf{%
\delta }}+\left( \mathbf{B}+\varepsilon \mathbf{C}\left( \mathbf{\delta }%
\right) \right) \mathbf{\delta }=0  \label{12}
\end{equation}%
where $\mathbf{A},\mathbf{B},\mathbf{C}\left( \mathbf{\delta }\right) ,%
\mathbf{D}$ are derived from the corresponding element matrices $%
A^{e},B^{e},C^{e}\left( \delta \right) ,D^{e}$, respectively and $\mathbf{%
\delta =}\left( \delta _{-1},...,\delta _{N+1}\right) ^{T}$ contain all
elements parameters.

The unknown parameters $\mathbf{\delta }$ are interpolated between two time
levels $n$ and $n+1$ with the Crank-Nicolson method%
\begin{equation*}
\begin{array}{cc}
\mathbf{\delta }=\dfrac{\delta ^{n+1}+\delta ^{n}}{2}, & \overset{\mathbf{%
\bullet }}{\mathbf{\delta }}=\dfrac{\delta ^{n+1}-\delta ^{n}}{\Delta t},%
\end{array}%
\end{equation*}%
we obtain iterative formula for the time parameters $\mathbf{\delta }^{n}$
\begin{equation}
\left[ \left( \mathbf{A}-\mu \mathbf{D}\right) \mathbf{+}\frac{\Delta t}{2}%
\left( \mathbf{B}+\varepsilon \mathbf{C}\left( \mathbf{\delta }^{n+1}\right)
\right) \right] \mathbf{\delta }^{n+1}=\left[ \left( \mathbf{A}-\mu \mathbf{D%
}\right) \mathbf{-}\frac{\Delta t}{2}\left( \mathbf{B}+\varepsilon \mathbf{C}%
\left( \mathbf{\delta }^{n}\right) \right) \right] \mathbf{\delta }^{n}
\label{13}
\end{equation}%
The set of equations consist of $\left( N+3\right) $ equations with $\left(
N+3\right) $ unknown parameters. Before starting the iteration procedure,
boundary conditions must be adapted into the system and initial vector\ $%
\mathbf{\delta }^{n}$ must also be determined.

We delete first and last equations from the system (\ref{13}) and eliminate
the terms $\delta _{-1}^{n+1}$ and $\delta _{N+1}^{n+1}$ from the system (%
\ref{13}) by using boundary conditions in (\ref{2}), which give the
following equations:

\begin{equation*}
\begin{array}{ll}
u\left( a,t\right) =m_{1}\delta _{-1}^{n}+\delta _{0}^{n}+m_{1}\delta
_{1}^{n}=\beta _{1}, & \text{ }u\left( b,t\right) =m_{1}\delta
_{N-1}^{n}+\delta _{N}^{n}+m_{1}\delta _{N+1}^{n}=\beta _{1}%
\end{array}%
\end{equation*}%
we obtain a septa-diagonal matrix with the dimension $\left( N+1\right)
\times \left( N+1\right) $. Then we can solve this matrix system through
Thomas algorithm. Since the system (\ref{13}) is an implicit system due to
the term $\mathbf{C}\left( \mathbf{\delta }^{n+1}\right) $, we have used the
following inner iteration:

\begin{equation}
(\mathbf{\delta }^{\ast }\mathbf{)}^{n+1}=\mathbf{\delta }^{n+1}+\dfrac{(%
\mathbf{\delta }^{n+1}-\mathbf{\delta }^{n})}{2}.  \label{14}
\end{equation}%
In this iteration, before moving the calculation of the next time step
approximation for time parameter, we calculate the new vectors $(\mathbf{%
\delta }^{\ast }\mathbf{)}^{n+1}$ using the formula (\ref{14}) from previous
vectors $\mathbf{\delta }^{n+1}$ finding form the system (\ref{13}) and then
repeat three times at all time steps.

To start evolution of the vector of initial parameters $\mathbf{\delta }^{0}$%
, it must be determined by using the initial condition and boundary
conditions:

\begin{equation}
\begin{tabular}{l}
$u_{0}^{\prime }(x_{0},0)=\dfrac{p\left( 1-c\right) }{2\left( phc-s\right) }%
\delta _{-1}+\dfrac{p\left( c-1\right) }{2\left( phc-s\right) }\delta _{1}$
\\
$u\left( x_{m},0\right) =\dfrac{s-ph}{2\left( phc-s\right) }\delta
_{m-1}+\delta _{m}+\dfrac{s-ph}{2\left( phc-s\right) }\delta _{m+1},$ $%
m=0,...,N$ \\
$u^{\prime }\left( x_{N},0\right) =\dfrac{p\left( 1-c\right) }{2\left(
phc-s\right) }\delta _{N-1}+\dfrac{p\left( c-1\right) }{2\left( phc-s\right)
}\delta _{N+1}$%
\end{tabular}
\label{15}
\end{equation}%
The solution of matrix equation (\ref{15}) with the dimensions $\left(
N+1\right) \times \left( N+1\right) $ is obtained by the way of Thomas
algorithm. Once $\mathbf{\delta }^{0}$ is determined, we can start the
iteration of the system to find the parameters $\mathbf{\delta }^{n}$ at
time $t^{n}=n\Delta t.$ Approximate solutions at the knots is found from the
Eq.(\ref{5}) and solution over the intervals $[x_{m},x_{m+1}]$ is determined
from the Eq.(\ref{15}).

\section{Test Problems}

We have carried out three test problems to demonstrate the given algorithm.
Accuracy of the method is measured by the error norm:%
\begin{equation*}
L_{\infty }=\left\Vert u^{\text{exact}}-u^{\text{numeric}}\right\Vert
_{\infty }=\max_{0\leq j\leq N}\left\vert u_{j}^{\text{exact}}-u_{j}^{\text{%
numeric}}\right\vert .
\end{equation*}%
The RLW equation satisfy the following conservation laws which are
corresponding to mass, momentum and energy \cite{olver}:%
\begin{equation*}
C_{1}=\int_{-\infty }^{+\infty }udx,C_{2}=\int_{-\infty }^{+\infty }\left(
u^{2}+\mu \left( u_{x}\right) ^{2}\right) dx,C_{3}=\int_{-\infty }^{+\infty
}\left( u^{3}+3u^{2}\right) dx.
\end{equation*}

In numerical calculations, the conservation laws are calculated by use of
the trapezoidal rule and the determination of $p$ in the exponential
B-spline is made by experimentally.

\subsection{Propagation of single solitary wave}

\begin{figure}[h]
\centering\includegraphics[scale=0.4]{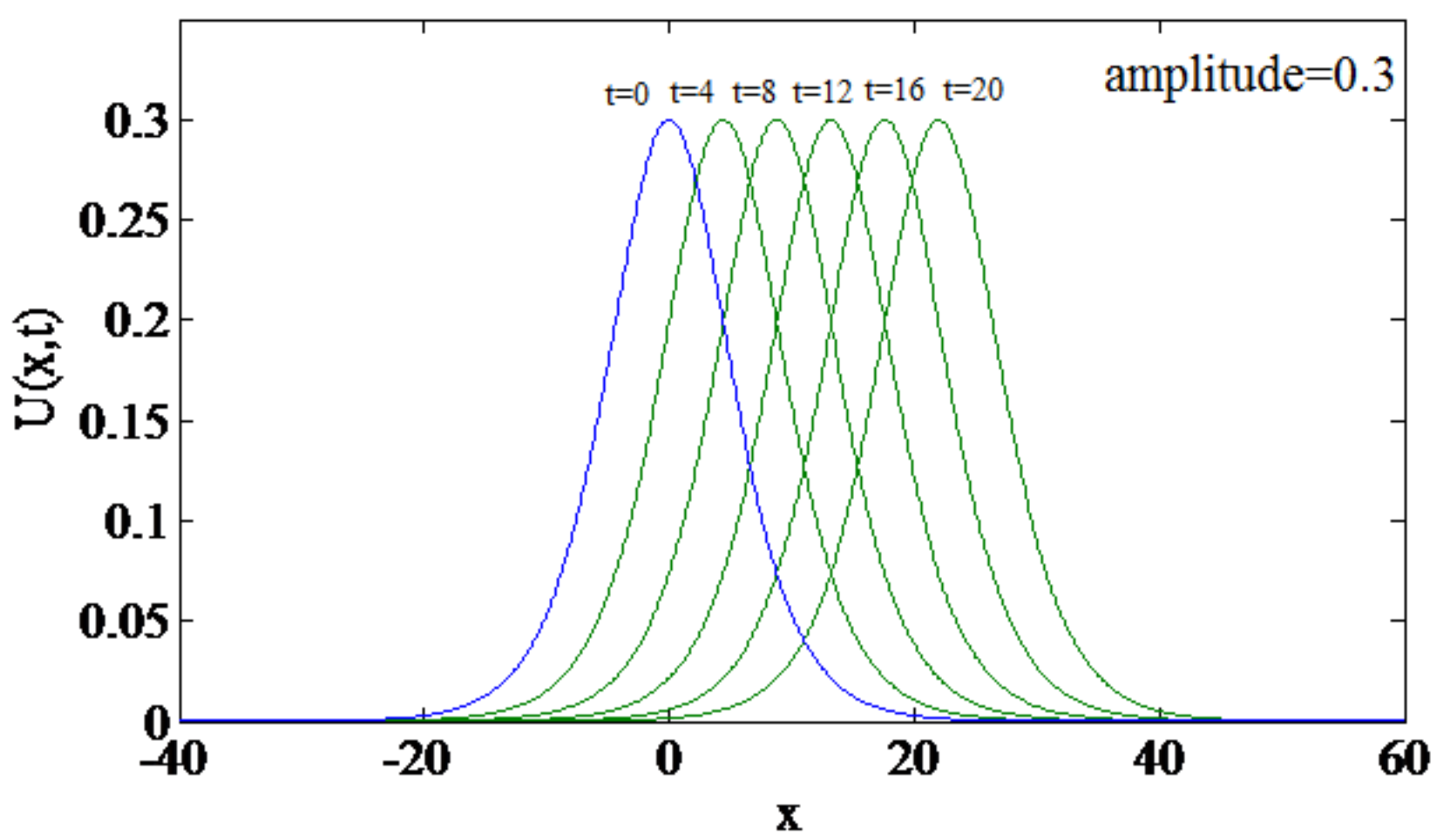}\includegraphics[scale=0.4]{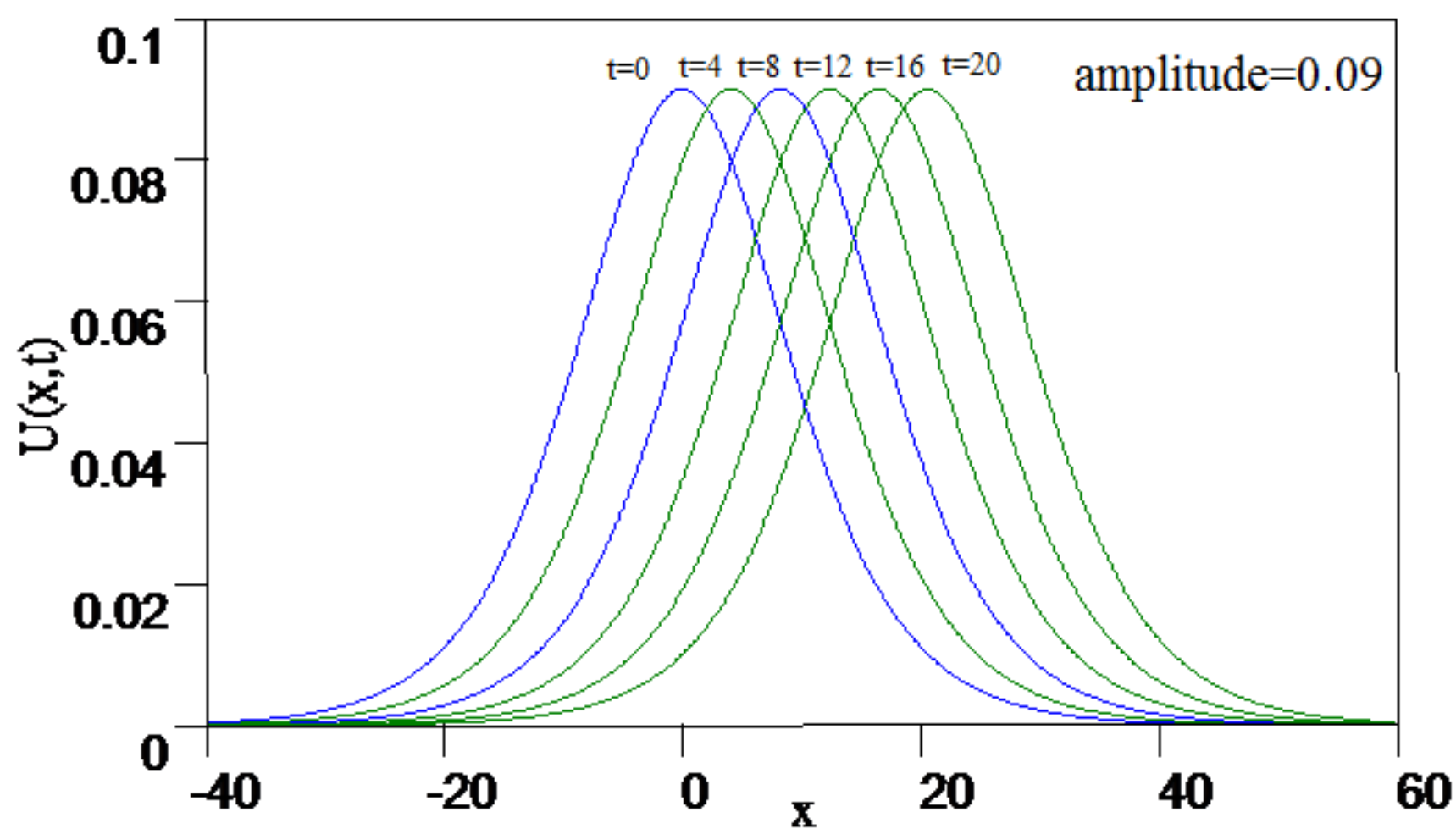}
\caption{{}{\protect\small \ Solitary wave solutions}}
\label{fig1-fig1b}
\end{figure}

The exact solution of RLW equation is given in \cite{peregrine} as follows:%
\begin{equation}
u(x,t)=\frac{3c}{\cosh ^{2}\left( k\left( x-x_{0}-\left( 1+\varepsilon
c\right) t\right) \right) }  \label{exact}
\end{equation}%
where $k=\dfrac{1}{2}\sqrt{\dfrac{\varepsilon c}{\mu \left( 1+\varepsilon
c\right) }}.$ This form of the solution is known as a single solitary wave
with the amplitude $3c$ and the velocity $1+\varepsilon c$. The initial
condition is obtained by taking $t=0$ in Eq.(\ref{exact})\textbf{. }We have
used boundary conditions $\beta _{1}=0,\beta _{2}=0.$The values of the
parameters seen in the above equations as%
\begin{equation*}
c=0.1\text{ and }0.03,\text{ }x_{0}=0,\text{ }\varepsilon =1,\text{ }\mu =1.
\end{equation*}%
\noindent With these parameters and the mentioned initial condition, the
solitary wave moves across the interval $-40\leq x\leq 60$ in time period $%
0\leq t\leq 20.$ Similar with some early studies, space step $h=0.125$ and\
time step $\Delta t=0.1$ are used in numerical calculations. In this test
problem, the $p=0.01262$ is determined by scanning the interval $[0,80]$
with the increment $0.1$ first, then according to the results scanning the
interval $[0,1]$ with the increment $0.00001$\textbf{. }The solution
profiles are illustrated in Figure \ref{fig1-fig1b} at selected times. It is
clear from this figure that the peak of the solitary wave remain kept during
the running time.

\begin{figure}[h]
\centering\includegraphics[scale=0.4]{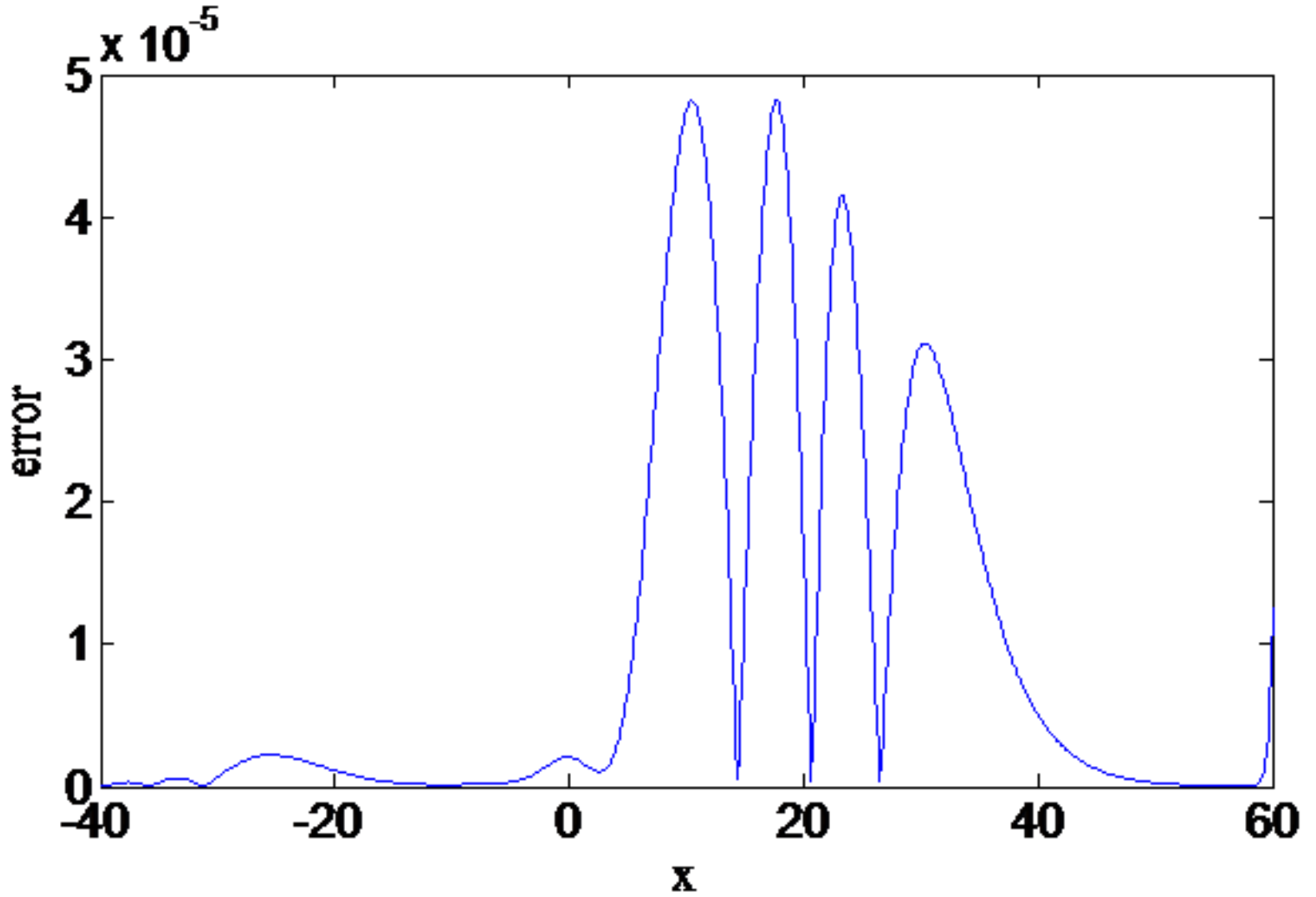}\includegraphics[scale=0.4]{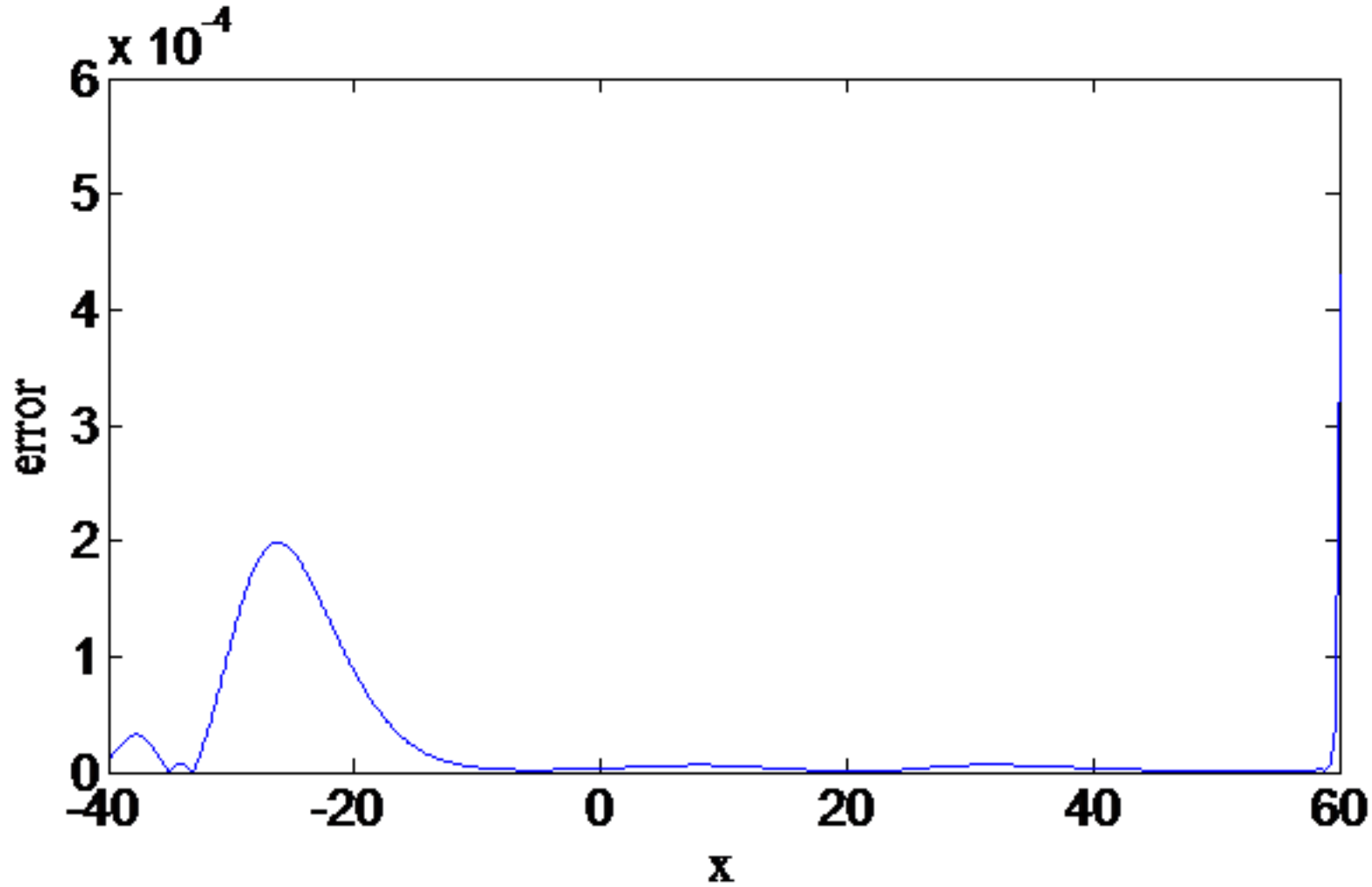}
\caption{{}{\protect\small amplitude= 0.3 and amplitude= 0.09}}
\label{fig2-fig3}
\end{figure}

The distribution of absolute error at $t=20$ for amplitude $0.3$ and $0.09$
is given in Figure \ref{fig2-fig3}, respectively. The maximum error for
\textbf{EBSGM} occurs at the right hand boundary seen in Figure \ref%
{fig2-fig3}. We believe that this error arises due to magnitude of the wave
and the physical boundary conditions to fit $u\left( a,0\right) \approx 0$
and $u\left( b,20\right) \approx 0$. If we extend the solution interval from
$\left[ -40,60\right] $ to $[-80,120],L_{\infty }$ error norm is seen to
reduce from $0$.$4315\times 10^{-3}$ to $0.63035\times 10^{-5}$ at time $%
t=20.$%

The absolute error norms and the values of the conservation invariants $%
C_{1},C_{2},C_{3}$ are recorded in Table 2 and 3 for different amplitudes.
To make a comparison with some early studies, the maximum errors with the
conservation invariants are presented in Table 4 and 5. According to this
tables, \textbf{EBSGM }is more accurate method than the some others. The
values of the conservation invariants $C_{1},C_{2},C_{3}$ at different times
remain fairly the same when compared with the analytical invariants $%
C_{1}=3.9799297,C_{2}=0.81046249,C_{3}=2.579007$ for amplitude $0.3.$ When
we take the amplitude as $0.09$, the value of the conservation invariant $%
C_{1}$ has some minor difference than the analytical value of it, whereas $%
C_{2}$ and $C_{3}$ are fairly same at different times.%
\begin{equation*}
\begin{tabular}{ccccc}
\multicolumn{5}{c}{{\small Table 2: Errors and Invariants for amplitude }$%
{\small 3c=0.3}${\small ,}} \\
\multicolumn{5}{c}{${\small h=0.125}${\small , }${\small \Delta t=0.1}$%
{\small , }${\small p=0.01262}$} \\ \hline\hline
{\small Time} & ${\small L}_{{\small \infty }}{\small \times }10^{5}$ & $%
C_{1}$ & $C_{2}$ & $C_{3}$ \\ \hline
\multicolumn{1}{l}{%
\begin{tabular}{r}
${\small 0}$ \\
${\small 4}$ \\
${\small 8}$ \\
\multicolumn{1}{l}{${\small 12}$} \\
\multicolumn{1}{l}{${\small 16}$} \\
\multicolumn{1}{l}{${\small 20}$}%
\end{tabular}%
} & \multicolumn{1}{l}{%
\begin{tabular}{l}
${\small 0.0000}$ \\
${\small 1.268}$ \\
${\small 2.445}$ \\
${\small 3.418}$ \\
${\small 4.197}$ \\
${\small 4.820}$%
\end{tabular}%
} &
\begin{tabular}{l}
${\small 3.9799267}$ \\
${\small 3.9799299}$ \\
${\small 3.9799281}$ \\
${\small 3.9799255}$ \\
${\small 3.9799168}$ \\
${\small 3.9798827}$%
\end{tabular}
&
\begin{tabular}{l}
${\small 0.8104625}$ \\
${\small 0.8104625}$ \\
${\small 0.8104625}$ \\
${\small 0.8104626}$ \\
${\small 0.8104626}$ \\
${\small 0.8104626}$%
\end{tabular}
&
\begin{tabular}{l}
${\small 2.5790074}$ \\
${\small 2.5790075}$ \\
${\small 2.5790076}$ \\
${\small 2.5790077}$ \\
${\small 2.5790078}$ \\
${\small 2.5790079}$%
\end{tabular}
\\ \hline
\end{tabular}%
\end{equation*}%
\begin{equation*}
\begin{tabular}{ccccc}
\multicolumn{5}{c}{{\small Table 3: Errors and Invariants for amplitude }$%
{\small 3c=0.09}${\small ,}} \\
\multicolumn{5}{c}{${\small h=0.125}${\small , }${\small \Delta t=0.1}$%
{\small , }${\small p=0.01262}$} \\ \hline\hline
{\small Time} & ${\small L}_{{\small \infty }}{\small \times 10}^{4}$ & $%
C_{1}$ & $C_{2}$ & $C_{3}$ \\ \hline
\begin{tabular}{r}
${\small 0}$ \\
${\small 4}$ \\
${\small 8}$ \\
\multicolumn{1}{l}{${\small 12}$} \\
\multicolumn{1}{l}{${\small 16}$} \\
\multicolumn{1}{l}{${\small 20}$}%
\end{tabular}
&
\begin{tabular}{l}
${\small 0.000}$ \\
${\small 2.302}$ \\
${\small 2.214}$ \\
${\small 2.129}$ \\
${\small 2.139}$ \\
${\small 4.315}$%
\end{tabular}
&
\begin{tabular}{l}
${\small 2.1070467}$ \\
${\small 2.1070975}$ \\
${\small 2.1068946}$ \\
${\small 2.1065454}$ \\
${\small 2.1053719}$ \\
${\small 2.1045885}$%
\end{tabular}
&
\begin{tabular}{l}
${\small 0.1273013}$ \\
${\small 0.1273011}$ \\
${\small 0.1273011}$ \\
${\small 0.1273011}$ \\
${\small 0.1273012}$ \\
${\small 0.1273012}$%
\end{tabular}
&
\begin{tabular}{l}
${\small 0.3888047}$ \\
${\small 0.3888041}$ \\
${\small 0.3888041}$ \\
${\small 0.3888040}$ \\
${\small 0.3888032}$ \\
${\small 0.3888023}$%
\end{tabular}
\\ \hline
\end{tabular}%
\end{equation*}%
\begin{equation*}
\begin{tabular}{lcccccc}
\multicolumn{7}{c}{\small Table 4:} \\
\multicolumn{7}{c}{{\small Errors and Invariants for amplitude }${\small %
3c=0.3}${\small ,} ${\small h=0.125}${\small ,} ${\small \Delta t=0.1}$%
{\small ,} ${\small p=0.01262}$} \\ \hline\hline
{\small Method} & {\small Time} & ${\small L}_{{\small \infty }}{\small %
\times 10}^{{\small 3}}$ & ${\small C}_{{\small 1}}$ & ${\small C}_{{\small 2%
}}$ & ${\small C}_{{\small 3}}$ &  \\ \hline
{\small Analytical} &  &  & ${\small 3.97995}$ & ${\small 0.8104624}$ & $%
{\small 2.579007}$ &  \\ \hline
\begin{tabular}{l}
{\small EBSGM} \\
{\small \cite{avilez}} \\
{\small \cite{dag3}} \\
{\small \cite{dag2}} \\
{\small \cite{saka}} \\
{\small \cite{saka2}} \\
{\small \cite{mei}} \\
{\small \cite{saka3}} \\
{\small \cite{dogan}}%
\end{tabular}
& \multicolumn{1}{l}{%
\begin{tabular}{r}
${\small 20}$ \\
${\small 20}$ \\
\multicolumn{1}{l}{${\small 20}$} \\
\multicolumn{1}{l}{${\small 20}$} \\
\multicolumn{1}{l}{${\small 20}$} \\
${\small 20}$ \\
${\small 20}$ \\
${\small 20}$ \\
${\small 20}$%
\end{tabular}%
} & \multicolumn{1}{l}{%
\begin{tabular}{l}
${\small 0.04820}$ \\
${\small 0.02643}$ \\
${\small 0.07337}$ \\
${\small 1.56640}$ \\
${\small 0.10299}$ \\
${\small 0.116}$ \\
${\small 0.91465}$ \\
${\small 0.07344}$ \\
${\small 0.073}$%
\end{tabular}%
} & \multicolumn{1}{l}{%
\begin{tabular}{l}
${\small 3.979883}$ \\
${\small 3.979909}$ \\
${\small 3.979883}$ \\
${\small 3.961597}$ \\
${\small 3.979858}$ \\
${\small 3.979883}$ \\
${\small 3.97972}$ \\
${\small 3.979888}$ \\
${\small 3.97989}$%
\end{tabular}%
} & \multicolumn{1}{l}{%
\begin{tabular}{l}
${\small 0.8104626}$ \\
${\small 0.8104625}$ \\
${\small 0.8104612}$ \\
${\small 0.804185}$ \\
${\small 0.8104596}$ \\
${\small 0.8102762}$ \\
${\small 0.81026}$ \\
${\small 0.8104622}$ \\
${\small 0.81046}$%
\end{tabular}%
} &
\begin{tabular}{l}
${\small 2.579008}$ \\
${\small 2.579007}$ \\
${\small 2.579003}$ \\
${\small 2.558292}$ \\
${\small 2.578999}$ \\
${\small 2.578393}$ \\
${\small 2.57873}$ \\
${\small 2.579006}$ \\
${\small 2.57901}$%
\end{tabular}
& \multicolumn{1}{l}{} \\ \hline
\end{tabular}%
\end{equation*}%
\begin{equation*}
\begin{tabular}{lccccc}
\multicolumn{6}{c}{\small Table 5:} \\
\multicolumn{6}{c}{{\small Errors and Invariants for amplitude }${\small %
3c=0.09}${\small ,} ${\small h=0.125}${\small ,} ${\small \Delta t=0.1}$%
{\small ,} ${\small p=0.01262}$} \\ \hline\hline
{\small Method} & {\small Time} & ${\small L}_{{\small \infty }}{\small %
\times 10}^{3}$ & ${\small C}_{1}$ & $C_{3}$ & $C_{3}$ \\ \hline
{\small Analytical} &  &  & \multicolumn{1}{l}{%
\begin{tabular}{l}
{\small 2.109407}%
\end{tabular}%
} & \multicolumn{1}{l}{%
\begin{tabular}{l}
{\small 0.12730171}%
\end{tabular}%
} & \multicolumn{1}{l}{%
\begin{tabular}{l}
{\small 0.3888059}%
\end{tabular}%
} \\ \hline
\begin{tabular}{l}
{\small EBSGM} \\
{\small \cite{dag2}} \\
{\small \cite{saka2}} \\
{\small \cite{mei}} \\
{\small \cite{saka3}} \\
{\small \cite{dogan}}%
\end{tabular}
& \multicolumn{1}{l}{%
\begin{tabular}{r}
{\small 20} \\
{\small 20} \\
\multicolumn{1}{l}{\small 20} \\
\multicolumn{1}{l}{\small 20} \\
\multicolumn{1}{l}{\small 20} \\
{\small 20}%
\end{tabular}%
} & \multicolumn{1}{l}{%
\begin{tabular}{l}
{\small 0.431512} \\
{\small 1.5506} \\
{\small 0.432} \\
{\small 0.439145} \\
{\small 0.19806} \\
{\small 0.199}%
\end{tabular}%
} & \multicolumn{1}{l}{%
\begin{tabular}{l}
{\small 2.104589} \\
{\small 2.128869} \\
{\small 2.104584} \\
{\small 2.10902} \\
{\small 2.104708} \\
{\small 2.10467}%
\end{tabular}%
} & \multicolumn{1}{l}{%
\begin{tabular}{l}
{\small 0.12730121} \\
{\small 0.127228} \\
{\small 0.12729366} \\
{\small 0.12730} \\
{\small 0.1273006} \\
{\small 0.12730}%
\end{tabular}%
} & \multicolumn{1}{l}{%
\begin{tabular}{l}
{\small 0.3888023} \\
{\small 0.388571} \\
{\small 0.3887776} \\
{\small 0.38880} \\
{\small 0.3888025} \\
{\small 0.38880}%
\end{tabular}%
} \\ \hline
\end{tabular}%
\end{equation*}%
\bigskip

\subsection{Interaction of two solitary waves}

In this section, we will study the interaction of two solitary waves having
different amplitudes and moving in the same direction. The initial condition
is%
\begin{eqnarray*}
u(x,0) &=&u_{1}+u_{2}, \\
u_{j} &=&\frac{3A_{j}}{\cosh ^{2}\left( k_{j}\left( x-\widetilde{x}%
_{j}\right) \right) },\text{ \ \ }A_{j}=\frac{4k_{j}^{2}}{1-4k_{j}^{2}},%
\text{ \ \ }j=1,2
\end{eqnarray*}%
where the following parameters are chosen to coincide values in the
literature:%
\begin{equation*}
\varepsilon =1,\text{ }\mu =1,\text{ }k_{1}=0.4,\text{ }k_{2}=0.3,\text{ }%
\widetilde{x}_{1}=15,\text{ }\widetilde{x}_{2}=35\text{ }\beta _{1}=0,\beta
_{2}=0.
\end{equation*}%
Initially, these parameters yields the solitary waves with the amplitudes $%
5.333375$ and $1.687502$ positioning around points $x=15$ and $x=35$
respectively. The computation is carried out up to time $t=30$ with time
step $\Delta t=0.1,$ $N=400$ over the finite interval $\left[ 0,120\right] $
and $p$ is selected as $1$ for \textbf{EBSGM.}

\begin{figure}[h]
\centering\includegraphics[scale=0.5]{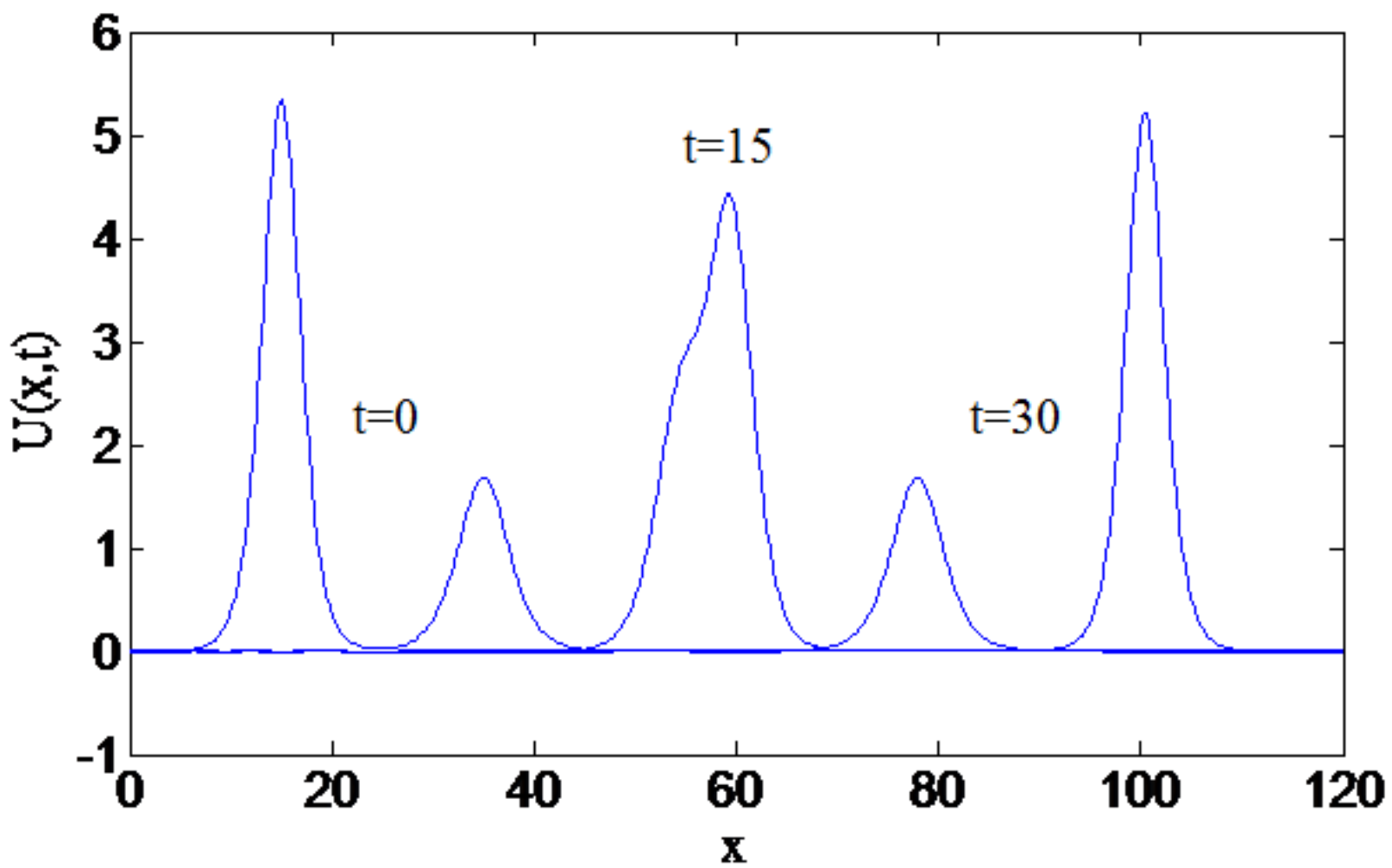}
\caption{{}{\protect\small Interaction of two positive solitary waves at }$%
{\protect\small t=0}${\protect\small , }${\protect\small 15}${\protect\small %
\ and }${\protect\small 30}${\protect\small .}}
\label{fig4}
\end{figure}

Numerical solutions of $u\left( x,t\right) $ at various times are depicted
in  Figure \ref{fig4},  and the initial solution has been propagated
rightward. It is seen from the  Figure \ref{fig4} that the solitary waves
are subjected to a collision about time $t=15$ and after the interaction
they propagate with their original amplitudes to the right seeing \ at time $%
t=30.$%

The conservation invariants are presented at some selected times in Table 6.
According to this table, during the interaction there are some changes at
the conservation constants $C_{2}$ and $C_{3}$ whereas the constant $C_{1}$
remain nearly same.

\begin{equation*}
\begin{tabular}{cccc}
\multicolumn{4}{c}{{\small Table 6: Invariants for }${\small h=0.3}${\small ,%
} ${\small \Delta t=0.1}${\small ,} ${\small p=1}${\small .}} \\ \hline\hline
{\small Time} & $C_{1}$ & $C_{2}$ & $C_{3}$ \\ \hline
\begin{tabular}{r}
{\small 0} \\
{\small 5} \\
{\small 10} \\
\multicolumn{1}{l}{\small 15} \\
\multicolumn{1}{l}{\small 20} \\
\multicolumn{1}{l}{\small 25} \\
{\small 30}%
\end{tabular}
& \multicolumn{1}{l}{%
\begin{tabular}{l}
{\small 37.916502} \\
{\small 37.861378} \\
{\small 37.810403} \\
{\small 37.790851} \\
{\small 37.752779} \\
{\small 37.701304} \\
\multicolumn{1}{r}{\small 37.650120}%
\end{tabular}%
} & \multicolumn{1}{l}{%
\begin{tabular}{l}
{\small 120.522247} \\
{\small 119.772286} \\
{\small 119.129874} \\
{\small 119.012803} \\
{\small 118.312087} \\
{\small 117.600957} \\
\multicolumn{1}{r}{\small 116.909863}%
\end{tabular}%
} & \multicolumn{1}{l}{%
\begin{tabular}{l}
{\small 744.081209} \\
{\small 737.819263} \\
{\small 732.319072} \\
{\small 730.931470} \\
{\small 725.698903} \\
{\small 719.868731} \\
\multicolumn{1}{r}{\small 714.174760}%
\end{tabular}%
} \\ \hline\hline
\end{tabular}%
\end{equation*}

\subsection{Wave generation}

An applied force like an introduction of fluid mass, an action of some
mechanical device, to a free surface, will induce waves. In this numerical
experiment, we take following boundary condition to generate waves with the
RLW equation.%
\begin{equation*}
u\left( a,t\right) =\beta _{1}=\left\{
\begin{array}{lcl}
U_{0}\dfrac{t}{\tau }, & \text{ \ \ \ } & 0\leq t\leq \tau , \\
U_{0}, &  & \tau <t<t_{0}-\tau , \\
U_{0}\dfrac{t_{0}-t}{\tau }, &  & t_{0}-\tau \leq t\leq t_{0}, \\
0, &  & \text{otherwise}.%
\end{array}%
\right.
\end{equation*}%
and $u(b,t)=\beta _{2}=0$ is studied to generate waves. This forced boundary
condition known as a wave maker at one end.

The parameters $U_{0}=2,$ $\Delta t=0.1,$ $h=0.4,$ $t_{0}=20,$ $\tau =0.3$
are chosen to make a comparison with earlier works over the region $0\leq
x\leq 260$. $p$ is selected as $1$ for \textbf{EBSGM. }During the run time
of the algorithm, five solitary waves are produced. Although first four
waves have reached amplitudes larger than forcing amplitudes, the last one
is less than that of the forcing one. When forcing is switched off, the last
wave has not enough time to evolve. Subsequently, no new wave are born. A
view of travelling solitary wave is presented at time $t=100$ in Figure \ref%
{fig5}. Amplitudes of solitary waves versus time are depicted in Figure \ref%
{fig6}. At various time, amplitudes of the solitary waves and the
conservation constants are demonstrated in Table 7 for \textbf{EBSGM. }In
addition to this, other amplitudes of the solitary waves which are reduced
the other studies are shown in Table 7 for time $t=100$. Our results are in
conformity with that of studies \cite{gardag2,chang,dogan}.
\begin{figure}[h]
\centering\includegraphics[scale=0.5]{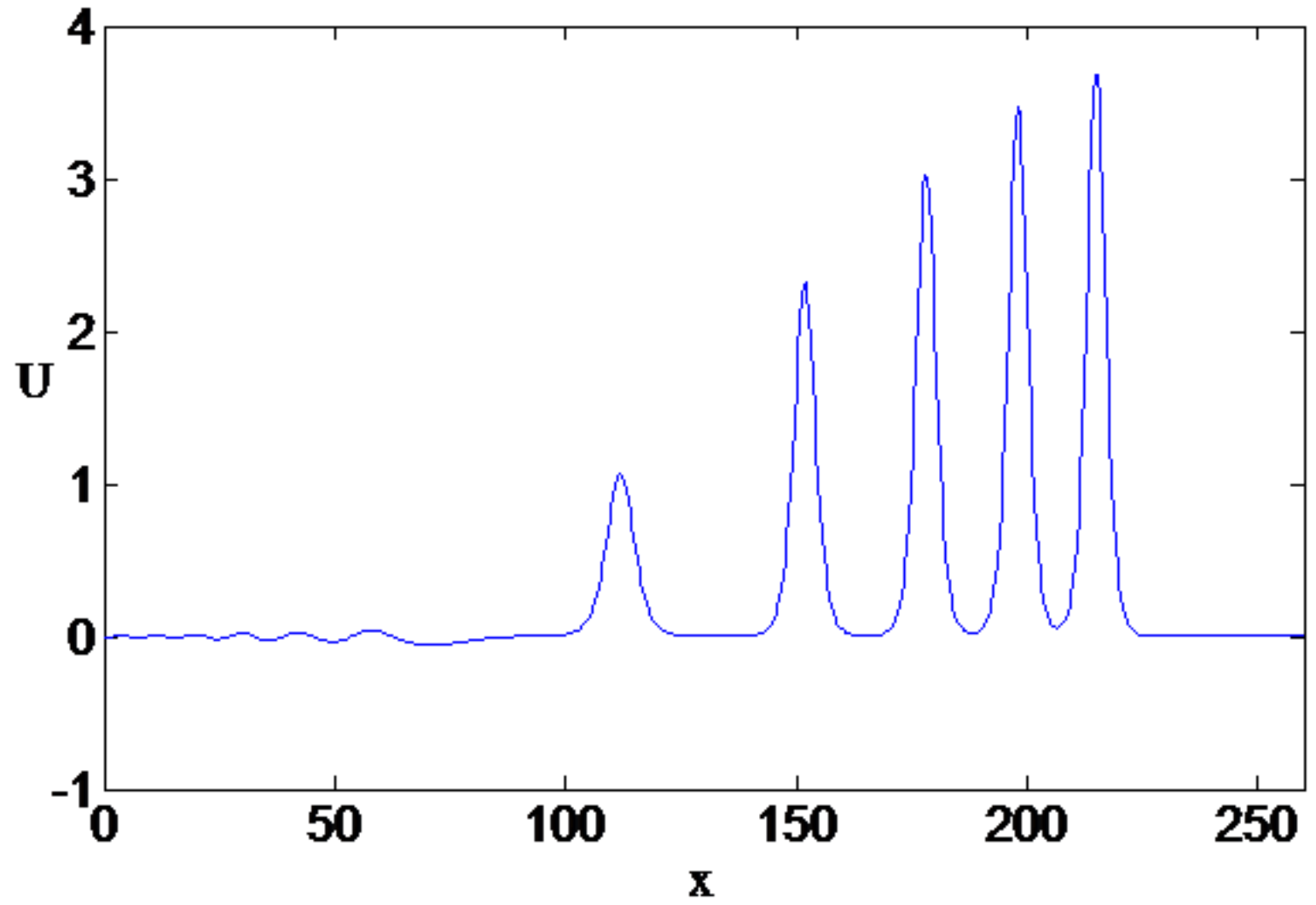}
\caption{{}{\protect\small Solitary wave produced by boundary forcing of
duration }${\protect\small t}_{{\protect\small 0}}{\protect\small =20}$
{\protect\small and amplitude} ${\protect\small U}_{{\protect\small 0}}%
{\protect\small =2}${\protect\small \ at time }${\protect\small t=100}$%
{\protect\small , }${\protect\small h=0.4}${\protect\small , }$%
{\protect\small \Delta t=0.1}${\protect\small , }${\protect\small p=1}$%
{\protect\small .}}
\label{fig5}
\end{figure}
\begin{figure}[h]
\centering\includegraphics[scale=0.5]{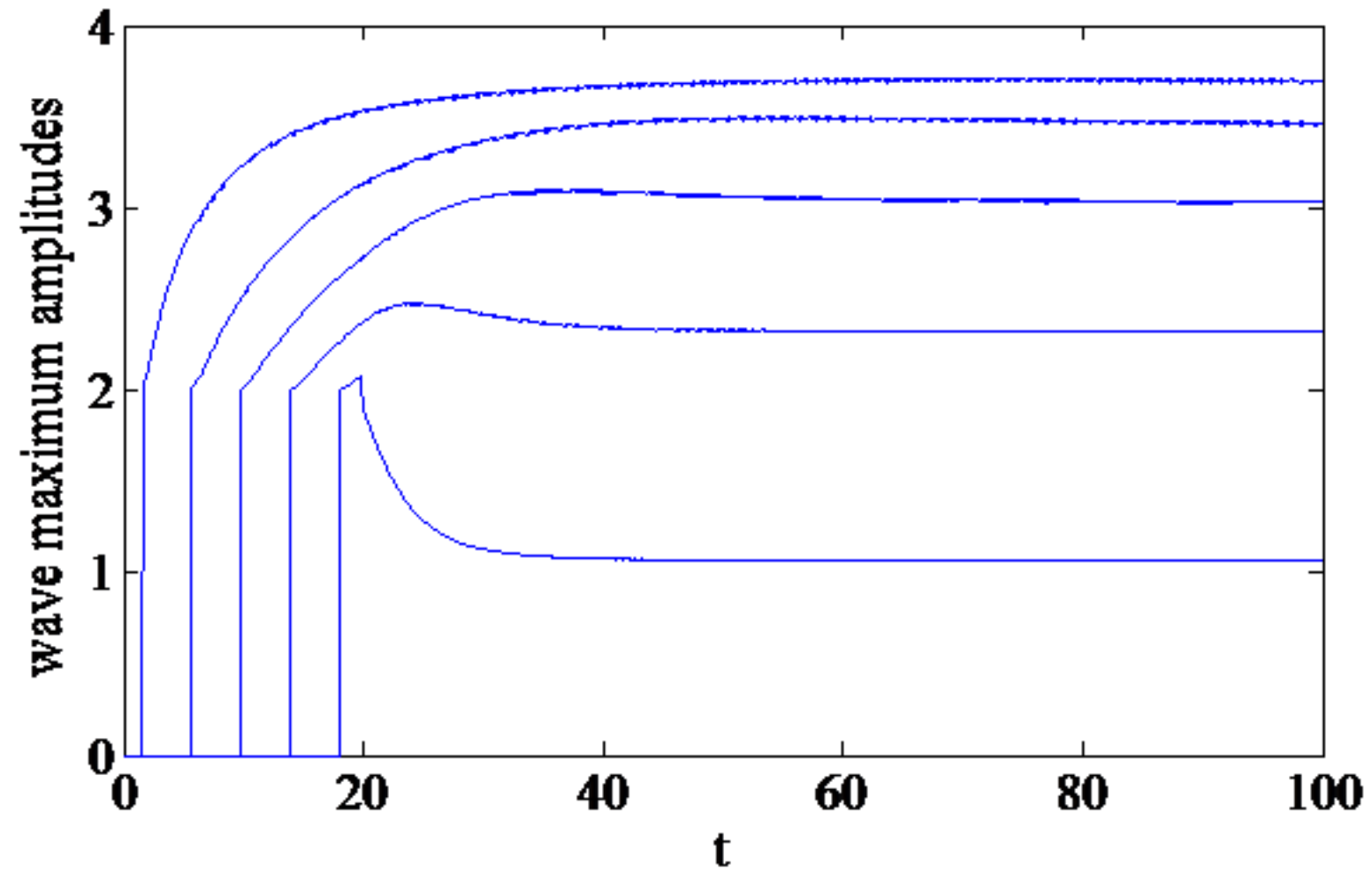}
\caption{{}{\protect\small Evolution of wave amplitudes.}}
\label{fig6}
\end{figure}
\begin{equation*}
\begin{tabular}{ccccccccc}
\multicolumn{9}{c}{$\text{{\small Table 7: Solitary wave amplitude with }}%
{\small U}_{{\small 0}}{\small =2}\text{{\small , }}{\small h=0.4}\text{%
{\small , }}{\small \Delta t=0.1}\text{{\small , }}{\small p=1}$} \\
\multicolumn{9}{c}{{\small and period of forcing }${\small 0\leq t}_{{\small %
0}}{\small \leq 20}$, ${\small 0\leq x\leq 260}${\small .}} \\ \hline\hline
{\small Time} & {\small 1} & {\small 2} & {\small 3} & {\small 4} & {\small 5%
} & $C_{1}$ & $C_{2}$ & $C_{3}$ \\ \hline
\multicolumn{1}{l}{%
\begin{tabular}{l}
${\small 2.5}$ \\
${\small 5}$ \\
${\small 7.5}$ \\
${\small 10}$ \\
${\small 15}$ \\
${\small 20}$ \\
${\small 40}$ \\
${\small 60}$ \\
${\small 80}$ \\
${\small 100}$ \\
${\small 100\cite{gardag2}}$ \\
${\small 100\cite{dogan}}$ \\
${\small 100\cite{chang}}$%
\end{tabular}%
} &
\begin{tabular}{l}
${\small 2.26}$ \\
${\small 2.78}$ \\
${\small 3.05}$ \\
${\small 3.23}$ \\
${\small 3.44}$ \\
${\small 3.53}$ \\
${\small 3.67}$ \\
${\small 3.70}$ \\
${\small 3.70}$ \\
${\small 3.68}$ \\
${\small 3.76}$ \\
${\small 3.77}$ \\
${\small 3.76}$%
\end{tabular}
&
\begin{tabular}{l}
\\
\\
${\small 2.24}$ \\
${\small 2.53}$ \\
${\small 2.91}$ \\
${\small 3.14}$ \\
${\small 3.46}$ \\
${\small 3.50}$ \\
${\small 3.48}$ \\
${\small 3.47}$ \\
${\small 3.52}$ \\
${\small 3.52}$ \\
${\small 3.51}$%
\end{tabular}
&
\begin{tabular}{l}
\\
\\
\\
${\small 2.02}$ \\
${\small 2.43}$ \\
${\small 2.73}$ \\
${\small 3.08}$ \\
${\small 3.05}$ \\
${\small 3.03}$ \\
${\small 3.03}$ \\
${\small 3.06}$ \\
${\small 3.08}$ \\
${\small 3.07}$%
\end{tabular}
&
\begin{tabular}{l}
\\
\\
\\
\\
${\small 2.06}$ \\
${\small 2.38}$ \\
${\small 2.35}$ \\
${\small 2.33}$ \\
${\small 2.33}$ \\
${\small 2.32}$ \\
${\small 2.33}$ \\
${\small 2.38}$ \\
${\small 2.31}$%
\end{tabular}
&
\begin{tabular}{l}
\\
\\
\\
\\
\\
${\small 1.90}$ \\
${\small 1.08}$ \\
${\small 1.07}$ \\
${\small 1.07}$ \\
${\small 1.07}$ \\
${\small 1.07}$ \\
${\small 1.18}$ \\
${\small 0.98}$%
\end{tabular}
& \multicolumn{1}{r}{%
\begin{tabular}{r}
${\small 8.9397}$ \\
\multicolumn{1}{l}{${\small 19.3958}$} \\
\multicolumn{1}{l}{${\small 29.3109}$} \\
\multicolumn{1}{l}{${\small 39.2124}$} \\
\multicolumn{1}{l}{${\small 59.2209}$} \\
\multicolumn{1}{l}{${\small 76.8502}$} \\
\multicolumn{1}{l}{${\small 78.5060}$} \\
\multicolumn{1}{l}{${\small 78.4275}$} \\
\multicolumn{1}{l}{${\small 78.3394}$} \\
\multicolumn{1}{l}{${\small 78.2589}$} \\
\\
\\
\end{tabular}%
} &
\begin{tabular}{r}
${\small 16.1942}$ \\
${\small 41.3363}$ \\
${\small 64.3064}$ \\
${\small 87.2152}$ \\
\multicolumn{1}{l}{${\small 133.8372}$} \\
\multicolumn{1}{l}{${\small 175.1059}$} \\
\multicolumn{1}{l}{${\small 174.3022}$} \\
\multicolumn{1}{l}{${\small 173.5305}$} \\
\multicolumn{1}{l}{${\small 172.7643}$} \\
\multicolumn{1}{l}{${\small 172.0094}$} \\
\\
\\
\end{tabular}
&
\begin{tabular}{r}
${\small 74.8986}$ \\
\multicolumn{1}{l}{${\small 206.5459}$} \\
\multicolumn{1}{l}{${\small 324.6697}$} \\
\multicolumn{1}{l}{${\small 442.2752}$} \\
\multicolumn{1}{l}{${\small 682.1665}$} \\
\multicolumn{1}{l}{${\small 882.7008}$} \\
\multicolumn{1}{l}{${\small 882.6289}$} \\
\multicolumn{1}{l}{${\small 877.6376}$} \\
\multicolumn{1}{l}{${\small 872.6609}$} \\
\multicolumn{1}{l}{${\small 867.7610}$} \\
\\
\\
\end{tabular}
\\ \hline
\end{tabular}%
\end{equation*}

\section{Conclusion}

In this paper, we investigate the utility of the exponential B-spline
algorithm for solving the RLW equation. The efficiency of the method is
tested on the propagation of the single solitary wave, the interaction of
two solitary waves and wave generation. To see the accuracy of the method, $%
L_{\infty }$ error norm and conservation quantities $C_{1},$ $C_{2}$ and $%
C_{3}$ are documented based on the obtained results. Exponential B-spline
based method gives accurate and reliable results for solving the RLW
equation. For the first test problem, \textbf{EBSGM }leads to more accurate
results than the collocation-based method but similar results with the some
Galerkin methods. In the second test problem, there is no exact solution
therefore simulation is shown graphically and the conservation quantities
are tabulated. Generation of waves by using variable boundary conditions at
left has been achieved and wave profiles and their amplitudes \ are
documented. In conclude, the numerical algorithm in which the exponential
B-spline functions are used, performs well compared with other existing
numerical methods for the solution of RLW equation.

\textbf{Acknowledgements}

The author, Melis Zor\c{s}ahin G\"{o}rg\"{u}l\"{u}, is grateful to The
Scientific and Technological Research Council of Turkey for granting
scholarship for PhD studies and all of the authors are grateful to The
Scientific and Technological Research Council of Turkey for financial
support for their project.

\end{document}